\documentclass[10pt,letter]{amsart}
\usepackage{amssymb,amsbsy,amsmath,amsfonts,amssymb}
\usepackage{latexsym,euscript,exscale}

\title{On the algebraic structure of the unitary group}
\author {\'Eric Ricard and Christian Rosendal}
\date {}
\newcommand {\A}{\mathbb A}
\newcommand {\B}{\mathbb B}
\newcommand {\F}{\mathbb F}

\newcommand {\N}{\mathbb N}

\newcommand {\U}{\mathbb U}

\newcommand{\W}{\mathbb W}

\renewcommand{\leq}{\ensuremath{\leqslant}}
\renewcommand{\geq}{\ensuremath{\geqslant}}
\newcommand{\ran}{\ensuremath{\mathop{\rm Im\,}}}
\renewcommand{\ker}{\ensuremath{\mathop{\rm Ker\,}}}
\renewcommand{\dim}{\ensuremath{\mathop{\rm dim\,}}}

\newcommand {\E}{\mathbb E}
\newcommand {\D}{\mathbb D}

\newcommand{\V}{\mathbb V}

\newcommand{\begr}{\!\upharpoonright}

\newcommand{\til}{\rightarrow}

\newcommand {\del}{ \; \big| \;}

\newcommand {\go} {\mathfrak}

\newcommand{\inv}{^{-1}}

\newcommand{\pff}{$\hfill  \Box$

\smallskip }

\newtheorem{thm}{Theorem}

\newtheorem{lemme}[thm]{Lemma}
\newtheorem{prop} [thm] {Proposition}

\begin{document}

\maketitle

\begin{abstract}
We consider the unitary group $\U$ of complex, separable, infinite-dimensional Hilbert space as a discrete group. It is proved that, whenever $\U$ acts by isometries on a metric space, every orbit is bounded. Equivalently, $\U$ is not the union of a countable chain of proper subgroups, and whenever $\E\subseteq \U$ generates $\U$, it does so by words of a fixed finite length. 
\end{abstract}

A property of uncountable groups that has recently been
studied by a number of authors is a strengthening of the
property of uncountable cofinality that originated in the work of
Jean-Pierre Serre on actions of groups on trees \cite{serre:trees}.
Here an uncountable group $G$ is said to have {\em uncountable
cofinality} if $G$ is not the union of a countable increasing chain
of proper subgroups. Serre proved this to be one of the three
conditions in his reformulation of when a group does not have fixed
point free actions without inversions on trees and it has
subsequently been confirmed for a great number of profinite groups
(see, e.g., Koppelberg and Tits \cite{koppelberg-tits}) and groups
of permutations of $\N$. The strengthening of this property, in
which we are interested, comes from considering the additional condition on $G$ that
whenever $E$ is a symmetric generating set for $G$ containing the
identity there is some finite $n$ such that $G=E^n$, i.e., any
element of $G$ can be written as a word of length $n$ using elements
of $E$. We denote this condition by {\em Cayley boundedness} and it
indeed corresponds to any Cayley graph for $G$ being of bounded
diameter in the word metric. George M. Bergman \cite{bergman:generating} originally
proved that the conjuction of uncountable cofinality and Cayley
boundedness holds for the infinite symmetric group $S_\infty$ and
subsequently this has been verified for a number of other
groups.  The purpose of this paper is to prove it for the unitary
group of complex separable infinite-dimensional Hilbert space
$\ell_2$. However, before we begin, we should mention some of the
equivalent formulations of these properties.

Consider a group $G$. It is fairly straightforward to see that $G$
has uncountable cofinality and is Cayley bounded if and only if it
satisfies the following property:
\begin{quotation}
Whenever $\W_0\subseteq \W_1\subseteq \W_2\ldots\subseteq
G=\bigcup_n\W_n$ is an increasing and exhaustive sequence of subsets
of $G$ there is some $n$ and some $k$ such that $\W^k_n=G$.
\end{quotation}
More interestingly is the fact that it is also equivalent to:
\begin{quotation}
Whenever $G$ acts by isometries on a metric space $(X,d)$ every orbit is bounded.
\end{quotation}
Or for another variation:
\begin{quotation}
Any left-invariant metric on $G$ is bounded.
\end{quotation}
There are several names for these equivalent properties in the
literature, but the two most common seem to be {\em Bergman
property} and {\em strong uncountable cofinality}.

The reformulation in terms of isometric actions discloses that we
are really dealing with a property of geometric group theory and
indeed it does have quite evident geometric implications for other
types of actions. For example, one can quite easily prove that any
group with this property must also have properties (FH) and (FA),
meaning that any action of the group by affine isometries on a real
Hilbert space has a fixed point and similarly, any action of the
group on a tree without inversions fixes a vertex. We should stress
the fact that despite that we will use a bit of topology in our
proof, the result is really about the unitary group as a {\em
discrete} group. Thus there are no topological assumptions being
made. We shall denote the unitary group of separable
infinite-dimensional complex Hilbert space $\ell_2$ by
$\U=\U(\ell_2)$. We recall that the inverse of a unitary operator is
its adjoint, $U\inv=U^*$. To avoid confusion, we will use the
blackboard bold font type $\A,\B,\ldots$. to denote subsets of $\U$
and standard fonts $A,B,\ldots, a,b,\ldots$ to denote individual
operators in $\go B(\ell_2)$. Thus $\A^*=\{A^*\del A\in \A\}$.

\begin{thm}\label{berg}
The unitary group $\U$ of separable, infinite-dimensional, complex
Hilbert space is Cayley bounded and has uncountable cofinality.
\end{thm}

The proof will proceed through a sequence of reductions showing that
the sets $\W_n$ of an exhaustive sequence of subsets are big in both
an algebraic and an analytic sense. We first need the following result.

\begin{thm}(Brown and Pearcy \cite{brown-pearcy}) Any $T\in
\U$ is a multiplicative commutator, i.e., $T=ABA^* B^* $ for
appropriate $A,B\in \U$.\end{thm}

\begin{prop}\label{bergman} Suppose $\W_0\subseteq \W_1\subseteq \ldots\subseteq
\U=\bigcup_n \W_n$. Then there is a decomposition $\ell_2=X\oplus
X^\perp$ into closed infinite-dimensional subspaces, such that for
some $l\in \N$,
$$
\D:=\{T\in \U\del T[X]=X\}=\U(X)\oplus \U(X^\perp)\subseteq
\W_l^{20}
$$
\end{prop}

\begin{proof} Suppose $\W_0\subseteq \W_1\subseteq \ldots\subseteq
\U=\bigcup_n\W_n$. By instead considering the sequence
$$
\W_0\cap \W_0^* \subseteq \W_1\cap \W_1^*\subseteq \ldots \subseteq
\U=\bigcup_n \W_n\cap \W_n^*
$$
we can suppose that each $\W_n$ is symmetric, $\W_n=\W^*_n$. Now,
write $\ell_2$ as the direct sum of infinitely many
infinite-dimensional closed subspaces $\ell_2=(\sum_n\oplus
X_n)_{\ell_2}$. Then for some $n$ and all $T\in \U$, if $T[X_n]=X_n$
then there is some $S\in \W_n$ with $T\begr X_n=S\begr X_n$. If not,
we would be able to find for each $n$ some unitary operator $T_n$ of
$X_n$ such that for all $S\in \W_n$, $T_n\neq S\begr X_n$. But then
the infinite direct sum $T=\oplus_m T_m$ is such that for all $n$
and $S\in \W_n$, $T\begr X_n=T_n\neq S\begr X_n$. In particular,
$T\notin \bigcup_n \W_n$, contradicting our supposition.

So suppose this holds for $n$. Find infinite-dimensional closed
subspaces $X,Y,Z\subseteq \ell_2$ such that $X=X_n$ and $X\oplus
Y\oplus Z=\ell_2$. Let also $N,M\in \U$ be such that $N[X]=X\oplus
Y$, $N[Y\oplus Z]=Z$, $M[X]=X\oplus Z$ and $M[Y\oplus Z]=Y$.

Clearly, if $A\in \U$ is such that $A[X\oplus Y]=X\oplus Y$, then
for some $S\in \W_n$, $S[X]=X$, we have
$$
A\begr X\oplus Y= NSN^* \begr X\oplus Y
$$
and similarly, if $B\in \U$, $B[X\oplus Z]=X\oplus Z$, then there is
$R\in \W_n$, $R[X]=X$, with
$$
B\begr X\oplus Z= MRM^* \begr X\oplus Z
$$
Now suppose that $T\in \U$, $T[X]=X$ and $T\begr X^\perp=T\begr
Y\oplus Z={\rm id} \begr X^\perp$. By the theorem of Brown and
Pearcy, there are unitary operators $A$ and $B$ on $X$ such that
$T\begr X=ABA^* B^*$. Extend now $A$ and $B$ to all of $\ell_2$ by
letting $A\begr X^\perp=B\begr X^\perp= {\rm id}\begr X^\perp$.

Moreover, find $S$ and $R$ in $\W_n$ as above, whence for $\hat
A=NSN^*$ and $\hat B=MRM^*$ we have
$$
\hat A\begr X\oplus Y=NSN^* \begr X\oplus Y= A\begr X\oplus Y
$$
$$
\hat B\begr X\oplus Z= MRM^* \begr X\oplus Z=B \begr X\oplus Z
$$
whence $\hat A\begr Y={\rm id}\begr Y$, $\hat B\begr Z={\rm id}\begr
Z$, $\hat A[Z]=Z$ and $\hat B[Y]=Y$. Thus
$$
\hat A\hat B\hat A^* \hat B^* \begr X=ABA^* B^* \begr X= T\begr X,
$$
$$
\hat A\hat B\hat A^* \hat B^* \begr Y=\hat B\hat B^*\begr Y={\rm
id}\begr Y= T\begr Y,
$$
and
$$
\hat A\hat B\hat A^* \hat B^* \begr Z=\hat A\hat A^*\begr Z={\rm
id}\begr Z= T\begr Z,
$$
Therefore,
\begin{align*}
T=\hat A\hat B\hat A^* \hat B^*&=(NSN^*)(MRM^*)(NS^* N^*) (MR^*
M^*)\\
&=NS(N^* M)R(M^* N)S^* (N^* M)R^* M^*\in \W_m^9
\end{align*}
provided that $m\geq n$ is large enough such that $N, M, N^* M\in
\W_m$. Notice however, that $N$ and $M$ do not depend on $T$, so
$\V=\{T\in \U\del T\begr X^\perp={\rm id }\begr X^\perp\}\subseteq
\W^9_m$. Find now a $K\in \U$ such that $K[X]=X^\perp$. Then
clearly, if $T\in \U$ satisfies $T[X]=X$, we have $T\in \V K\V K^*$.
Let now $l\geq n$ be sufficiently big that $K\in \W_l$. Then we have
$\D\subseteq \W_l^9\W_l\W^9_l\W_l=\W_l^{20}$.
\end{proof}

The {\em strong operator topology} on $\U$ is the topology of
pointwise convergence on $\ell_2$, i.e., $U_i\til U$ if for all
$x\in \ell_2$, $U_ix\til Ux$. In this topology $\U$ becomes a {\em
Polish} space, i.e., a separable space whose topology is induced by
a complete metric. Actually $\U$ is a $G_\delta$ in $\go B(\ell_2)$
under this topology. A subset $A$ of a Polish space $X$ is said to
be {\em analytic} if it is the image of another Polish space by a
continuous function. Analytic sets have the Baire property, meaning
that they differ from an open set by a meagre set (see, e.g.,
Kechris \cite{kechris:classical} for the basics of descriptive set
theory). The following result was proved in \cite{rosendal:OB} as a
byproduct of other computations, but for the readers convenience we
include a simple proof here.

\begin{prop}\label{finitude}
Assume that $\F$ is a symmetric subset of $\U$, closed in the strong
operator topology, and that $U_0,U_1,U_2,\ldots$ is a sequence of
unitary operators such that $\U$ is generated as a group by
$\{U_n\}$ over $\F$. Then there is some finite $n$ such that
$$(\F\cup\{1, U_0, U_0^*,\ldots,U_n,U_n^*\})^n=\U.$$
\end{prop}

\begin{proof}
Define $\W_n=(\F\cup\{1, U_0, U_0^*,\ldots,U_n,U_n^*\})^n$ and
notice that $(\W_n)$ is an increasing, exhaustive sequence of
symmetric, analytic subsets of $\U$. Since the sets have the Baire
property there is some $k$ such that $\W_k$ is comeagre in an open
set and hence by Pettis' Theorem (see Kechris
\cite{kechris:classical}), $\W_k\W_k^*=\W_k^2$ contains an open
neighbourhood of the identity in the strong operator topology. This
implies that we can find some finite dimensional space $X\subseteq
\ell_2$ such that if $U\begr X=1$, then $U\in \W_k^2$. Now find a
unitary operator $V$ such that $Y:=V[X]\subseteq X^\perp$. Then
$V\W_k^2V^*$ contains all unitaries $U$ such that $U\begr Y=1$.
Suppose now that $U$ is an arbitrary unitary operator and find some
finite dimensional space $Z\subseteq (X\oplus Y)^\perp$ such that
$U[X]\subseteq X\oplus Y\oplus Z$. Find now some $W_0\in \W^2_k$ such that
$W_0[Y]\subseteq (X\oplus Y\oplus Z)^\perp$, while $W_0\begr X\oplus Z=I$. Then
$$
W_0U[X]\subseteq W_0[X\oplus Y\oplus Z]\subseteq X\oplus W_0[Y]\oplus Z\subseteq Y^\perp.
$$
There is therefore some $W_1\in V\W^2_kV^*$ such that $W_1W_0U\begr
X=I$. Thus we get $W_1W_0U\in \W_k^2$ and therefore
$$
U\in W_0^*
W_1^* \W^2_k\subseteq \W_k^2V\W^2_kV^*\W_k^2
$$
and
$$
\U=\W_k^2V\W^2_kV^*\W_k^2.
$$
Find now some sufficiently big $m>k$
such that $V,V^* \in \W_m$, then
$$\U=\W_k^2V\W^2_kV^*\W_k^2=\W_k^2\W_m\W_k^2\W_m\W_k^2=\W_m^8=\W_{8m}.$$
So $n=8m$ works.
\end{proof}

We now fix some increasing, exhaustive sequence $(\W_n)$ of subsets
of $\U$, and note that by considering instead $(\W_n\cap\W_n^*)$ we
can assume they are symmetric. By Proposition \ref{bergman}, we also
fix a decomposition $\ell_2=X\oplus X^\perp$ into closed
infinite-dimensional subspaces such that
$\D=\U(X)\oplus \U(X^\perp)\subseteq \W_l^{20}$ for some $l$.

\begin{prop}\label{uetd}
$\U$ is finitely generated over $\D$.
\end{prop}

We notice that this will indeed be enough to prove Theorem
\ref{berg}. For combining Propositions \ref{finitude} and
\ref{uetd}, we get that $\U$ has bounded length with respect to $\D$
and the family $U_1,U_1^*,...,U_8,U_8^*$ that we will prove
generates $\U$ over $\D$, in other words, there is an integer $k$
such that any unitary $U$ in $\U$ can be written as a product
$U=V_1...V_k$ where $V_i$ either belongs to $\D$ or is one of $U_i$
or $U_i^*$.

We also note also that, though we have not so far been able to
calculate it, the $k$ in the definition of the Bergman property does
not depend on the sequence $(\W_i)$. It is known that this
uniformity does not hold in all groups with the Bergman property,
as, e.g., it fails in the full group of a countable measure
preserving equivalence relation (see B.D. Miller
\cite{miller:dissertation}).

From now on, we will use identifications of both $X$ and $X^\bot$
with $\ell_2=\ell_2(\N)$ with its canonical basis and see
endomorphisms of $X\oplus X^\bot$ as two by two matrices with
entries in $\go B(\ell_2)$.

A partial isometry is a map $u$ on $\ell_2$ so that  $u^*u=p$ and
$uu^*=q$ are orthogonal projections. We say that $u$ is a partial
isometry with initial space $\ran p$ and final space $\ran q$ (or
from $\ran p$ to $\ran q$), where $\ran T$ denotes the closure of
the range of an operator $T$, and notice that, in this case, $u$ is
actually an isometric bijection between $\ran p$ and $\ran q$. We
recall also the polar decomposition of an operator: Every operator
$T$ can be decomposed as $T=u |T|$, where $|T|=(T^*T)^{1/2}$ is
positive, self-adjoint, $\ker |T|=\ker T$, and $u$ is a partial
isometry from $\ran |T|=(\ker T)^\bot$ to $\ran T$.

Let $S$ be the shift operator on $\ell_2$, and $L$ be an isometry
from $\ell_2$ to an infinite codimensional subspace ($L(e_i)=e_{2i}$
for instance). For any isometry $u$, $p_u$ will be the projection
$1-uu^*$. Let $\W$ be the subgroup of $\U$ generated by $\D$ and the
matrices
$$
\left[\begin{array}{cc}
 u & p_u \\
 0 & u^*
\end{array}\right],\, \frac 1 {\sqrt 2}
\left[\begin{array}{cc}
 1  & u^* \\
 u & -uu^*+\sqrt 2 p_u
\end{array}\right],\,  \frac 1 {\sqrt 2}
\left[\begin{array}{cc}
 1  & u^{*2} \\
 u^2 & -u^2u^{*2}+\sqrt 2 p_{u^2}
\end{array}\right], \,
\left[\begin{array}{cc}
 0 & 1 \\
 1 &  0
\end{array}\right],$$
where $u$ is either $1$, $S$ or $L$.

$\W$ is generated by $\D$ and eight unitaries, that we denote by
$U_1,..., U_8$. We will show that $\W=\U$.

\begin{lemme}\label{isom}
Let $\left[\begin{array}{cc}
 A  \\
 B
\end{array}\right]$ be an isometry from $\ell_2$ to $\ell_2\oplus \ell_2$, then
there is a partial isometry $v$ from $\ran 1-A^*A$ to $\ran B$
so that $B=v(1-A^*A)^{1/2}$.
\end{lemme}

\begin{proof}
Since $\left[\begin{array}{cc}
 A  \\
 B
\end{array}\right]$ is an isometry, we get that
$A^*A+B^*B=1$. So we conclude that $|B|=(1-A^*A)^{1/2}$, the last part
follows from the polar decomposition and the fact that for any positive
operator $T$,
$\ran T=\ran T^{1/2}$.
\end{proof}

\begin{lemme}\label{part}
$\W$ contains all matrices of the form
$$\left[\begin{array}{cc}
 u & p_u \\
 0 & u^*
\end{array}\right],\qquad \frac 1 {\sqrt 2}
\left[\begin{array}{cc}
 1  & u^* \\
 u & -uu^*+\sqrt 2 p_u
\end{array}\right],$$
where $u$ is an isometry.
\end{lemme}

\begin{proof}
First, remark that up to multiplications by unitaries any isometry
$u$ is determined only by the dimension of $\ran p_u$. So
multiplying by elements in $\D$, we only need to prove the lemma for
$u=S^k$, $u=L$ and $u=1$. The only non-trivial cases are for
$u=S^k$. For the first case, we notice that
$$\left[\begin{array}{cc}
 S & p_S \\
 0 & S^*
\end{array}\right]^k= \left[\begin{array}{cc}
 S^k & A \\
 0 & S^{*k}
\end{array}\right].$$
Thus $\left[\begin{array}{cc}
 A  \\
 S^{*k}
\end{array}\right]$
and $\left[\begin{array}{cc}
 S^{*k}  \\
 A^*
\end{array}\right]$
are isometries from $\ell_2$ to $\ell_2\oplus\ell_2$. Hence, as
$S^{*k}$ is an isometry on $(\ran
p_{S^{k}})^\perp=[e_{k+1},e_{k+2},\ldots]$, we have $A\begr (\ran
p_{S^{k}})^\perp=A^*\begr(\ran p_{S^{k}})^\perp=0$, while, as
$S^{*k}\begr \ran p_{S^{k}}=0$, $A$ and $A^*$ must be isometries on
$\ran p_{S^{k}}$. Therefore, $A$ is a partial isometry with initial
and final space $\ran p_{S^{k}}$. Letting $U=(1-p_{S^k}) + A$, we
see that $U$ is unitary and that
$$\left[\begin{array}{cc}
 S & p_S \\
 0 & S^*
\end{array}\right]^k\left[\begin{array}{cc}
 1 & 0 \\
 0 & U^*
\end{array}\right]= \left[\begin{array}{cc}
 S^k & p_{S^k} \\
 0 & S^{*k}
\end{array}\right].$$

For the other case, we notice that by definition $\frac 1 {\sqrt 2}
\left[\begin{array}{cc}
 1  & S^{*k} \\
 S^k & -S^kS^{*k}+\sqrt 2 p_{S^k}
\end{array}\right]$
is in $\W$  for $k=1$ and $k=2$. Now
\begin{displaymath}\begin{split}
&\left[\begin{array}{cc}
 S & p_S \\
 0 & S^*
\end{array}\right]^* \left[\begin{array}{cc}
 1  & S^{*k} \\
 S^k & -S^kS^{*k}+\sqrt 2 p_{S^k}
\end{array}\right]\left[\begin{array}{cc}
 S & p_S \\
 0 & S^*
\end{array}\right]
=\left[\begin{array}{cc}
 1  & S^{*(k+2)} \\
 S^{k+2} & X
\end{array}\right],\end{split}\end{displaymath}
where $X=p_S+p_SS^{*(k+1)}+S^{k+1}p_s-S^{k+1}S^{*(k+1)}+\sqrt 2 Sp_SS^*$.

Let $V=-S^{k+2}S^{*(k+2)}+\sqrt 2 p_{S^{k+2}}$ and let $u$ be the
unitary operator given by $u(e_1)=\frac{e_1+e_{k+2}}{\sqrt2}$,
$u(e_{k+2})=\frac {e_1-e_{k+2}}{\sqrt2}$, while $u(e_i)=e_i$ for all
$i\neq 1,k+2$. One can check, by, e.g., considering the action on
vectors, that $X=uV$ and, as $u$ is the identity on
$[e_1,e_{k+2}]^\perp$, $u^*S^{k+2}=S^{k+2}$. Thus
$$\left[\begin{array}{cc}
 1  & S^{*(k+2)} \\
 S^{k+2} & V
\end{array}\right]=\left[\begin{array}{cc}
 1 & 0 \\
 0 & u^*
\end{array}\right]\left[\begin{array}{cc}
 1  & S^{*(k+2)} \\
 S^{k+2} & X
\end{array}\right],$$
and the lemma follows by induction.
\end{proof}

The matrices appearing the preceding lemma are actually special
cases of the more general form of a unitary that will appear in the
following.

\medskip

\noindent{\it Proof of Proposition \ref{uetd}:} We have to prove
that any unitary belongs to $\W$. The strategy is to start with an
arbitrary unitary viewed as a 2 by 2 matrix
$U=\left[\begin{array}{cc}
 T & M \\
 N & K
\end{array}\right]$ and to multiply it by elements in $\W$ to get simpler
forms.

First, we can assume that $T$ is positive. To see this, notice that
by replacing $U$ by $U^*$ we can assume that $\dim \ker T\leq \dim
\ker T^*$. We have $T=u|T|$ for some partial isometry $u$ from $\ran
|T|=(\ker T)^\bot$ to $\ran T=(\ker T^*)^\bot$, and our hypothesis
$\dim \ker T\leq \dim \ker T^*$ implies that $u$ can be extended to
an isometry $\tilde u$ so that $T=\tilde u |T|$. Then
$$\left[\begin{array}{cc}
 \tilde u & p_{\tilde u} \\
 0 & \tilde{u}^*
\end{array}\right]^* U=\left[\begin{array}{cc}
|T| & * \\
 * & *
\end{array}\right],$$
which proves our claim.

Applying lemma \ref{isom} twice, we only need to prove that a
unitary of the form
$$U=\left[\begin{array}{cc}
T & (1-T^2)^{1/2}v^* \\
u(1-T^2)^{1/2} & X
\end{array}\right]$$
is in $\W$, where $u$ and $v$ are partial isometries with initial
space $\ran (1-T^2)^{1/2}=(\ker (1-T^2))^\bot=(\ker(1-T))^\perp$ (as $T$ is 
positive). So
$u^*u=v^*v$ is the projection onto $\ran (1-T^2)^{1/2}$. Let
$uu^*=1-p$ and $vv^*=1-q$ and note that $u^*p=pu=v^*q=qv=0$.

The fact that $U$ is a unitary means that $U^*U=UU^*=1$, and writing
down the computation on matrices, this implies that
\begin{eqnarray*}
&(i) & u(1-T^2)u^* + XX^*=1, \\
&(ii) & v(1-T^2)v^* + X^*X=1, \\
&(iii) & u(1-T^2)^{1/2} T +X v(1-T^2)^{1/2}=0, \\
&(iv) &  v(1-T^2)^{1/2} T +X^* u(1-T^2)^{1/2}=0.
\end{eqnarray*}
Now $u^*u=v^*v$ is the projection onto $\ran (1-T^2)^{1/2}$, so as
$T$ and $(1-T^2)^{1/2}$ commute, we have by $(iii)$
$$
uT(1-T^2)^{1/2}+Xv(1-T^2)^{1/2}=0,
$$
whence $-uT$ and $Xv$ agree on $\ran(1-T^2)^{1/2}=\ran v^*$. So
$$
-uTv^*=Xvv^*=X(1-q).
$$
Similarly, using $(iv)$ we have
$$
-vTu^*=X^*uu^*=X^*(1-p).
$$
Thus $-uTv^*=X(1-q)=(1-p)X$.

Using $(ii)$ we have
$|X|=(X^*X)^{1/2}=(1-v(1-T^2)v^*)^{1/2}=(q+vTv^*)^{1/2}$. Now, 
$q$ is the projection onto $(\ran v)^\perp$, so $qvT^2v^*=0=vT^2v^*q$,
and, using that $m^{1/2}$ is the operator limit of polynomials in
$m$, we therefore have (as $v^*v$ commutes with $T$)
$$ (q+vT^2v^*)^{1/2}=q^{1/2}+(vT^2v^*)^{1/2}=q+vTv^*$$

Hence, by the polar decomposition, 
there is a partial isometry $\alpha$ with initial space
$\ran(q+vTv^*)^{1/2}$ so that $X=\alpha(q+vT^2v^*)^{1/2}$.
The initial space of $\alpha$ is $\ran q+vTv^*=\ran
q + \ran vTv^*$ and therefore $\alpha^*\alpha\geq q$. Now,
$$X=Xq+X(1-q)=\alpha(q+vTv^*)q-uTv^*= \alpha q - uTv^*.$$

Similarly, using $(i)$, we see that $X=(p+uTu^*)\beta^*$, where
$\beta$ is a partial isometry with initial space $\ran p+uTu^*=\ran
p + \ran uTu^*$, whence $\beta^*\beta\geq p$, and
$$X=pX+(1-p)X= p\beta^* - uTv^*.$$
So $p\beta^*=\alpha q=pX=Xq=pXq$.

Using $\alpha^*\alpha\geq q$, we see that $\alpha q=p\beta^*$ is an
isometry with initial space $\ran q$, and, as $\beta^*\beta\geq p$,
$\alpha q=p\beta^*$ has final space $\ran p$.

To summarise,
$$U=\left[\begin{array}{cc}
T & (1-T^2)^{1/2}v^* \\
u(1-T^2)^{1/2} & -uTv^* + \delta
\end{array}\right],$$
where $u$ is a partial isometry from $(\ker 1-T)^\bot$ to $\ran
1-p$, $v$ is  a partial isometry from $(\ker 1-T)^\bot$ to $\ran
1-q$ and $\delta=\alpha q=p\beta^*$ is a partial isometry from $\ran
p$ to $\ran q$.

Consequently, one can find a unitary $s$ so that $v=su$  and that
coincides with $\delta$ on $\ran p$. So
$$
U'=U\left[
   \begin{array}{cc}
     1 & 0 \\
     0 & s \\
   \end{array}
 \right]=
\left[\begin{array}{cc}
T & (1-T^2)^{1/2}u^* \\
u(1-T^2)^{1/2} & -uTu^* + p
\end{array}\right].$$
In particular, $U'$ is self-adjoint, and note that $\ker 1-uTu^*=\{0\}$.

The next step is to show that we can suppose that $u$ is an
isometry.  To do this, it suffices to write $T'=T-q'$
where $q'$ is the projection onto $\ker 1-T$. So then $\ker
1-T'=\{0\}$. Then $T'q'=0$, $q'T'=q'-q'T=q'-(Tq')^*=q'-q'=0$, and $q'u^*=0$, so
$(1-T^2)^{1/2}u^*=(1-T^{\prime 2})^{1/2}u^*$ and 
$$U'=\left[\begin{array}{cc}
T'+q' & (1-T'^2)^{1/2}u^* \\
u(1-T'^2)^{1/2} & -uT'u^* + p
\end{array}\right]$$
where $u$ is a partial isometry from $\ran 1-q'$ onto $\ran 1-p$ and $\ker 1-T'=\{0\}$.

The previous formula is now symmetric in terms of $T$ and $uTu^*$ as 
$\ker 1-T'=\ker 1-uT'u^*=\{0\}$. We mean that  $T=T'+q'$ and 
$-uTu^*+p=-uT'u^*+p$ play exactly the same role in the above formula. So, 
using a conjugation with $\left[\begin{array}{cc}
 0 & 1 \\
 1 &  0
\end{array}\right]$ to reverse the roles of $p$ and $q'$, 
we therefore can assume that $\dim \ran p\geq \dim \ran q'$.
 Using other conjugations with elements in $\D$, we can  assume that
$q'\leq p$ and  $u$ is a partial isometry from $1-q'$ to $1-p$.
 
With these choices, it is
then easy to find an isometry $w$ extending $u$, that is $(w-u)(1-q')=0$, so
that replacing $T'+q'$ by $T'$, we can suppose that $U'$ is of the form
$$
\left[\begin{array}{cc}
T' & (1-T'^2)^{1/2}u^* \\
u(1-T'^2)^{1/2} & -uT'u^* + p
\end{array}\right]$$
were $u$ is an isometry from $\ell_2$ to $\ran 1-p$. Note that the matrices
of Lemma \ref{part} are of this form.

$T'$ is a self-adjoint positive contraction on $\ell_2$, so using the functional
calculus, define $A=\exp(i \arccos T')$ and $B=\exp(-i \arccos T')$. These
are two unitaries with 
$$
A+B=2 \cos(\arccos T')=2 T' \quad \textrm{ and} \quad 
A-B=2i\sin(\arccos T')=2i(1-T'^2)^{1/2}.
$$ Notice also that
$uBu^*-p$ is  unitary, actually $p=p_u$ with the notations of
Lemma \ref{part}. Finally,
\begin{eqnarray*} & & \frac 1 {\sqrt 2} \left[\begin{array}{cc}
1 & u^*  \\
u & -uu^* + \sqrt 2 p
\end{array}\right].\left[\begin{array}{cc}
A & 0  \\
0 & uBu^*-p
\end{array}\right].\frac 1 {\sqrt 2} \left[\begin{array}{cc}
1 & u^*  \\
u & -uu^* + \sqrt 2 p
\end{array}\right]\\&=&\left[\begin{array}{cc}
T' & i(1-T'^2)^{1/2}u^* \\
iu(1-T'^2)^{1/2} & uT'u^* - p
\end{array}\right]
\\ & =& \left[\begin{array}{cc}
1 & 0  \\
0 & -i
\end{array}\right].\left[\begin{array}{cc}
T' & (1-T'^2)^{1/2}u^* \\
u(1-T'^2)^{1/2} & -uT'u^* + p
\end{array}\right].\left[\begin{array}{cc}
1 & 0  \\
0 & -i
\end{array}\right]
\end{eqnarray*}
So $U\in \W$. \pff

There should be a proof of the bounded length of $\U$ that avoids
the Baire category argument.

\begin{flushleft}
{\em Address of \'E. Ricard:}\\
Laboratoire de Math\'ematiques,\\
Universit\'e de Franche-Comt\'e,\\
16 Route de Gray,\\
25030 Besan{\c c}on Cedex\\
France

\texttt{eric.ricard@math.univ-fcomte.fr}
\end{flushleft}

\begin{flushleft}
{\em Address of C. Rosendal:}\\
Department of mathematics,\\
University of Illinois at Urbana-Champaign,\\
273 Altgeld Hall, MC 382,\\
1409 W. Green Street,\\
Urbana, IL 61801,\\
USA.\\
\texttt{rosendal@math.uiuc.edu}
\end{flushleft}

\bibliography{RosendalGeneral}

\begin{thebibliography}{Kec95}

\bibitem[Ber]{bergman:generating}
George~M. Bergman.
\newblock Generating infinite symmetric groups.
\newblock To appear in Bull. London Math. Soc.

\bibitem[BP66]{brown-pearcy}
Arlen Brown and Carl Pearcy.
\newblock Multiplicative commutators of operators.
\newblock {\em Canad. J. Math.}, 18:737--749, 1966.

\bibitem[Kec95]{kechris:classical}
Alexander~S. Kechris.
\newblock {\em Classical descriptive set theory}, volume 156 of {\em Graduate
  Texts in Mathematics}.
\newblock Springer-Verlag, New York, 1995.

\bibitem[KT74]{koppelberg-tits}
Sabine Koppelberg and Jacques Tits.
\newblock Une propri\'et\'e des produits directs infinis de groupes finis
  isomorphes.
\newblock {\em C. R. Acad. Sci. Paris S\'er. A}, 279:583--585, 1974.

\bibitem[Mil04]{miller:dissertation}
Benjamin~David Miller.
\newblock {\em Full groups, classification, and equivalence relations}.
\newblock PhD thesis, UC Berkeley, 2004.

\bibitem[Ros]{rosendal:OB}
Christian Rosendal.
\newblock A topological version of the {B}ergman property.
\newblock Preprint.

\bibitem[Ser03]{serre:trees}
Jean-Pierre Serre.
\newblock {\em Trees}.
\newblock Springer Monographs in Mathematics. Springer-Verlag, Berlin, 2003.
\newblock Translated from the French original by John Stillwell, Corrected 2nd
  printing of the 1980 English translation.

\end{thebibliography}

\bibliographystyle{alpha}

\end{document}